\documentclass[11pt,reqno]{amsart}
\usepackage{amsmath,amssymb,amsthm,amsfonts,mathrsfs,latexsym,times,color,bm}
\usepackage{amsmath,amsfonts}
\usepackage{amsthm}
\usepackage{graphicx}
\usepackage{color}
\usepackage{multicol}
\usepackage{float}

\usepackage{bm}

\usepackage{hyperref}

\newtheorem{thm}{Theorem}[section]
\newtheorem{lem}{Lemma}[section]

\newtheorem{remark}{Remark}[section]

\newcommand{\kp}{\kappa}
\newcommand{\p}{\partial}

\newcommand{\R}{\mathbb R}
\newcommand{\bA}{\mathbf A}
\newcommand{\ba}{\mathbf a}

\newcommand{\bH}{\mathbf H}
\newcommand{\bG}{\mathbf G}

\def\charf {\mbox{{\text 1}\kern-.24em {\text l}}}

\numberwithin{equation}{section}

\begin{document}
	\title[Chern-Simons gauged $O(3)$ sigma equations]{Local well-posedness for Chern-Simons gauged $O(3)$ sigma equations under the Lorenz gauge}
	
	
	\author{Guanghui Jin}
	\address{Department of Mathematics, \newline Yanbian University, Yanji, 133002, China}
	\email{jinguanghui@ybu.edu.cn}
	
	\author{Huali Zhang}
	\address{Department of Mathematics, \newline Hunan University, Lushan South Road in Yuelu District, Changsha, 410882, China.}
	
	\email{hualizhang@hnu.edu.cn}

	\date{\today}


\begin{abstract}
In this paper, we study the Cauchy problem for the Chern-Simons gauged $O(3)$ sigma model under the Lorenz gauge condition. We prove the local well-posedness of solutions if the initial matter field and gauge field satisfy $(\bm{\phi}_0, \bA_0) \in H^s(\R^2)\times H^{s-\frac12}(\R^2)$, $s>1$, where the critical regularity for $\bm{\phi}_0$ is $s_c=1$. Our proof is based on identifying null forms within the system and utilizing bilinear estimates in wave-Sobolev space.
\end{abstract}
\maketitle



\section{Introduction}
\setcounter{equation}{0}

\ \\
The Chern-Simons guaged $O(3)$ sigma equation is a condensed matter systems describing the quantum Hall effects, high-temperature superconductors, and planar magnetism. The Lagrangian for the Chern-Simons gauged $O(3)$ sigma model on $\R^{1+2}$, proposed in \cite{GG, KLL}, is given by\footnote{Greek indices, such as $\mu,\,\nu,\,\rho$ will refer to all indices $0,\, 1,\, 2,$ whereas latin indices, such
	as $i,\, j,\, k,$ will refer only to the spatial indices $1,\, 2,$ unless otherwise specified. Repeated indices are summed over. $\epsilon^{\alpha\beta\gamma}$ is the totally-antisymmetric tensor with $\epsilon^{012}=1$. }
\begin{align}\label{lag}
	\mathcal{L}= \frac{1}{2}D_\mu\bm{\phi}\cdot D^\mu\bm{\phi}+\frac{\kappa}{4}\epsilon^{\mu\nu\rho} A_{\mu} F_{\nu \rho} -\frac{1}{2\kappa^2}(1+n_3\cdot\bm{\phi})(1-n_3\cdot \bm{\phi})^3,
\end{align}
where $\bm{\phi}$ is a three component vector $\bm{\phi}=({\phi}_1,{\phi}_2,{\phi}_3)^{\mathrm{T}}$ with the unit norm $\langle \bm{\phi}, \bm{\phi} \rangle=1$, where $\langle \cdot, \cdot \rangle$ stands for the inner product\footnote{For two vectors $a$ and $b$, the notation $a \cdot b$ also represents the inner product. We also use $|a|=\sqrt{\langle a, a \rangle}.$}. The functions $A_\mu:\R^{1+2}\rightarrow\R$ are the gauge field. Set $\bA=(A_0,A_1,A_2)^{\mathrm{T}}$. The gauge covariant derivative is defined by $D_\mu\bm{\phi}=\p_\mu\bm{\phi}+A_\mu (n_3\times\bm{\phi})$ and the Maxwell field is given by
$F_{\mu\nu}=\p_\mu A_\nu-\p_\nu A_\mu$, where $n_3=(0,0,1)^{\mathrm{T}}$. The constant $\kappa>0$ is a Chern-Simons coupling constant. The velocity of light and the Plank's constant in units of $\frac{1}{2\pi}$ are taken to be unity. In this case the coefficient of the CS term $(\kappa)$ has dimension of the inverse mass.

The factor $\frac{1}{2\kappa^2}$ in front of the potential term in \eqref{lag} is chosen so as to have a Bogomo'nyi bound. So we assume
\begin{equation}\label{c}
	0<\frac{1}{2\kappa^2} \leq m,
\end{equation}
where $m$ is a positive constant.

The Euler-Lagrange equation associated with \eqref{lag}, which is 
\begin{align}\tag{CSS}
\begin{aligned}\label{E-L}
D_\mu D^\mu \bm{\phi}+\bm{\phi}(D_\mu  \bm{\phi} \cdot D_\mu \bm{\phi})&=-\frac{1}{\kp^2}\left(\bm{\phi}(n_3\cdot\bm{\phi})-n_3(\bm{\bm{\phi}}\cdot\bm{\phi})\right)(1-n_3\cdot\bm{\phi})^2(1+2n_3\cdot\bm{\phi}),\\
\kp F_{\nu\lambda}&=-\epsilon_{\mu\nu\lambda}(\bm {\phi}\times D^\mu\bm{\phi})\cdot n_3.
\end{aligned}
\end{align}
\subsection{Some fundamental properties of CSS}

\ \\
The system \eqref{E-L} has conservation of energy[see Appendix for detail],
\begin{equation}\label{energy}
	E(\bm{\phi}, A)=\frac12 \int_{\R^2}\left( |D_\mu\bm{\phi}|^2 +\frac{1}{\kappa^2}(1+ n_3\cdot\bm{\phi})(1-n_3\cdot\bm{\phi})^3\right)dx.
\end{equation}

This system is invariant under the scaling
\[
\phi(t,x)\rightarrow\phi^{\lambda}:=\phi(\lambda t,\lambda x),\quad A_{\mu}\rightarrow A_{\mu}^{\lambda}:=\lambda A_{\mu}(\lambda t,\lambda x).
\]


By this scaling, we can calculate
\begin{equation*}
	\begin{split}
		\|\phi^{\lambda}\|_{\dot{H}^{s_c}(\mathbb{R}^2)} =& \lambda^{s_c-1} \|\phi\|_{\dot{H}^{s_c}(\mathbb{R}^2)},
		\\
		\|A_{\mu}^{\lambda}\|_{\dot{H}^{s'_c}(\mathbb{R}^2)} =& \lambda^{s_c'} \|A_{\mu}\|_{\dot{H}^{s'_c}(\mathbb{R}^2)}.
	\end{split}
\end{equation*}
Thus the critical Sobolev exponent for $\bm{\phi}$ and $\bA$ are $s_c=1$ and $s_c'=0$, respectively. From \eqref{energy}, we say that the Cauchy problem of \eqref{E-L} is energy critical. Therefore, we should consider the local well-posedness of \eqref{E-L} when the initial regularity is greater than the critical exponent. 

Another important property of the (CSS) is the gauge invariance. The system is invariant under the following gauge transformations:
$$\bm{\phi}=(z,{\phi}_3)\rightarrow(ze^{i\chi},{\phi}_3), \quad A_\mu\rightarrow A_\mu-\p_\mu\chi,$$
where $\chi$ is a real-valued smooth function on $\R^{1+2}$ and we use the notation $ z={\phi}_1+i {\phi}_2$. Therefore a solution of the system \eqref{E-L} is formed by a class of gauge equivalent pairs $(\bm{\phi},A_\mu).$ 
In this paper, we fix the gauge by imposing the Lorenz gauge condition
$\partial_\mu A^\mu=0$. Under this gauge, we can rewrite \eqref{E-L} as follows,
\begin{align}
	D_\mu D^\mu\bm{\phi}+\bm{\phi}(D_\mu\bm{\phi}\cdot D_\mu\bm{\phi})&=-\frac{1}{\kp^2}\left(\bm{\phi}(n_3\cdot\bm{\phi})-n_3(\bm{\phi}\cdot\bm{\phi})\right)(1-n_3\cdot\bm{\phi})^2(1+2n_3\cdot\bm{\phi}) ,\label{f0}\\
	\kappa F_{01}&=-\langle n_3 \times\bm{\phi}, D_2\bm{\phi} \rangle,\label{f1}\\
	\kappa F_{12}&=-\langle n_3 \times\bm{\phi}, D_0\bm{\phi} \rangle,\label{f2}\\
	\kappa F_{02}&=\langle n_3 \times\bm{\phi}, D_1\bm{\phi} \rangle,\label{f3}\\
	\partial_\mu A^\mu&=0\label{lorenz},
\end{align}
with the initial data
\begin{align}\label{initial}
	A_{\mu}(0,\cdot)=a_{\mu}, \quad \partial_t A_\mu(0,\cdot)=a'_\mu, \quad \bm{\phi}(0,\cdot)=\bm{\phi}_0,\quad \p_t\bm{\phi}(0,\cdot)=\bm{\phi}_1.
\end{align}

We can verify that the constraint $|\bm{\phi}|^2=1$ is preserved as follows. Define $\rho=|\bm{\phi}|^2-1$. Due to \eqref{f0}, $\rho$ satisfies
\begin{equation}\label{ro}
	\begin{split}
			&\left(  \partial_{\mu}  \partial^{\mu}+ 2 \langle D_{\mu}\bm{\phi}, \,D^{\mu} \bm{\phi}\rangle + 2{\phi}_3 (1-{\phi}_3)^2(1+2{\phi}_3) \right)\rho
		\\
		= & 2 \langle D_{\mu}\bm{\phi}, \,D^{\mu} \bm{\phi}\rangle+2 \langle \bm{\phi}, \,D_{\mu}D^\mu \bm{\phi}\rangle +2 |\bm{\phi}|^2\langle D_{\mu}\bm{\phi}, \,D^{\mu} \bm{\phi}\rangle-2 \langle D_{\mu}\bm{\phi}, \,D^{\mu} \bm{\phi}\rangle\\
		& \,\,\,\, +2{\phi}_3 |\bm{\phi}|^2(1-{\phi}_3)^2(1+2{\phi}_3)-2{\phi}_3(1-{\phi}_3)^2(1+2{\phi}_3)
		\\
		=& 0.
	\end{split}
\end{equation}
Seeing \eqref{ro}, it's a linear Klein-Gordon equation for $\rho$ with external potential $2 \langle D_{\mu}\bm{\phi}, \,D_{\mu} \bm{\phi}\rangle+2{\phi}_3(1-{\phi}_3)^2(1+2{\phi}_3)$. Given the initial data $\rho|_{t=0}=|\bm{\phi}_0|^2-1=0$ and
$\partial_t\rho|_{t=0}=2\langle\bm{\phi}_0 ,\, \bm{\phi}_1 \rangle=0$, we conclude that $\rho\equiv0$ over the interval $[0,\,T]\times \R^2$, thereby preserving the constraint $|\bm{\phi}|^2=1$.

Futhermore, the equation \eqref{f2} is propagetes along the system \eqref{E-L}. To verify this, we calculate
\begin{align}\label{constraint}
	\begin{aligned}
		&\partial_t \big( \kappa (\p_1 A_2-\p_2 A_1)+\langle n_3 \times\bm{\phi}, D_0\bm{\phi} \rangle \big)
		\\
		= & \kappa\partial_1F_{02}-\kappa\partial_2F_{01}+\langle D_0(n_3 \times\bm{\phi}), D_0\bm{\phi} \rangle+\langle n_3\times \bm \phi, D_0D_0\bm{\phi} \rangle \\
		=& \langle D_\mu( n_3 \times\bm{\phi}), D^\mu\bm{\phi} \rangle+\langle ( n_3 \times\bm{\phi}), D_\mu D^\mu\bm{\phi} \rangle \\
		=& \langle ( n_3 \times\bm{\phi}), D_\mu D^\mu\bm{\phi} \rangle
		\\
		=& 0.
	\end{aligned}
\end{align}
Thus, \eqref{constraint} implies that constraint \eqref{f2} is automatically satisfied at $t\geq0$ if the initial data satisfy
\begin{equation}\label{constraint1}
	\kappa(\p_1a_2-\p_2a_1)+\langle n_3\times\bm{\phi},\,\bm{\phi}_1+a_0(n_3 \times\bm{\phi}_0)\rangle=0.
	\end{equation}

To prove our main theorem, we will reformulate the system \eqref{f0}--\eqref{lorenz} as a wave equation incorporating null forms. Following \cite{FK}, we define the basic null forms: \begin{align*}
	 Q_0(u,,v)&=\partial_t u \partial_t v - \nabla u \cdot \nabla v,\\
	 Q_{ij}(u,,v)&=\partial_i u \partial_j v - \partial_j u \partial_i v,\\
	 Q_{0j}(u,,v)&=\partial_t u \partial_j v - \partial_j u \partial_t v,
	  \end{align*} 
	  where $\partial_j$ denotes spatial derivatives and $\nabla$ represents the spatial gradient.

Using these null forms, we can rewrite the system  \eqref{f0}--\eqref{lorenz} as follows: 
\begin{align} 
	&\Box \bm{\phi} = -\bm{\phi} \left( Q_0(\bm{\phi}, \bm{\phi}) + 2 A_{\mu} \p^{\mu} \bm{\phi} \cdot (n_3 \times \bm{\phi}) + A_{\mu} A^{\mu} |n_3 \times \bm{\phi}|^2 \right) - 2 A^{\mu} (n_3 \times \p_{\mu} \bm{\phi}) \label{main1} \\
	& \qquad \quad  - A_\mu A^\mu n_3 \times (n_3 \times \bm{\phi}) - \frac{1}{\kappa^2} (\phi_3\bm{\phi}  - n_3 )(1 -  {\phi}_3)^2 (1 + 2 {\phi}_3 ), \nonumber \\
	&\Box A_{\mu} = \frac{1}{\kappa} \epsilon_{\mu\nu\rho}\left( Q^{\nu \rho}(\bm{\phi}, \bm{\phi} \times n_3) + A^{\rho} \p^{\nu} \left( (\bm{\phi} \times (n_3 \times \bm{\phi})) \cdot n_3 \right) \right). \label{main2}
	 \end{align}

Additionally, if $\bA=0$ and $\frac{1}{\kappa^2}\rightarrow 0$, then the system \eqref{main1}--\eqref{main2} reduces to the wave map equation. We now present some historical results related to this context.

\subsection{Historical results}
The Cauchy problem for the wave map equation has been a central topic in the study of nonlinear hyperbolic partial differential equations due to its deep connections to geometry and mathematical physics. For wave maps, the energy critical sapce is $\dot H^{\frac n2}(\mathbb R^n)$. One important topic is about the low regularity local well-posedness, which started from the celebrated work due to Klainerman and Machedon. For $n=3$, Klainerman and Machedon \cite{KM1} proved the bilinear estimates of null forms. Later, Klainerman and Machedon \cite{KM2} also established the local well-posedness if $s>\frac32$ for the null-form quadratic wave equations. For $n\geq 2$, the local well-posedness for $s>\frac{n}{2}$ of wave maps was established by Klainerman and Selberg \cite{KS}. If $n=1$, Keel and Tao \cite{KT} showed that the wave maps is locally well-posed for $s>\frac12$. We note there are some other important low regularity results on wave maps due to Foschi-Klainerman \cite{FK}, Grigoryan-Nahmod \cite{GN}, Selberg \cite{Selberg}, Tao \cite{Tao,Tao1}, Zhou \cite{zhou1}. For global solutions of wave maps, we refer the readers to the significant progress by Shatah-Struwe \cite{SS}, Tataru \cite{ta,tata}, Krieger-Schlag-Tataru \cite{KST} and so on.

We should also mention another related model, i.e., Chern-Simons-Higgs (CSH) equation,
\begin{equation*}
	\begin{cases}
		& F_{\mu \nu}=  \frac{2}{\kappa} \epsilon_{\mu\nu\lambda} \mathrm{Im}( \bar{{\phi}} \mathbf{D}^\lambda {\phi}),\qquad (t,x)\in \mathbb{R}^{1+2},
		\\
		& \mathbf{D}^\mu  \mathbf{D}^\mu {\phi}=0,
		\\
		& ({\phi},\bA)|_{t=0}=({\phi}_0,\bA_0),
	\end{cases}
\end{equation*}
where the scalar field $\phi$ is a complex-valued function,  the 1-form $A_\mu$ represents a connection, and $F_{\mu \nu}:=\partial_\mu A_\nu- \partial_\nu A_\mu$ is the curvature 2-form, $\mathbf{D}_\mu:=\partial_\mu- \mathrm{i}A_\mu$ stands for the associated covariant derivative, and $\kappa>0$ is the Chern-Simons coupling constant. For the initial fields ${\phi}_0$ and $\bA_0$, the critical regularity are $s_c=\frac12$ and $s'_c=0$ respectively. Under the Coulomb gauge condition, the global existence of solutions to CSH for $(\bm{\phi}_0,\bA_0) \in H^2 \times H^{1}$ was shown by Chae and Choe \cite{CC}. Under the Lorenz
gauge, Selberg and Tesfahun \cite{ST} proved a low regularity local result for $(\bm{\phi}_0,\bA_0) \in H^{\frac{7}{8}+}\times H^{\frac{3}{8}+}$, and improved the global result of \cite{CC} with the data $(\bm{\phi}_0,\bA_0) \in H^1\times \dot{H}^{\frac12}$. Huh and Oh (\cite{HO}) further established the local well-posedness of CSH for $(\bm{\phi}_0,\bA_0) \in H^{\frac{3}{4}+}\times H^{\frac{1}{4}+}$.

Let us go back to Chern-Simons $O(3)$ gauged sigma models \eqref{E-L}. For the 1-dimensional case, Huh and Jin \cite{HJ} proved the local well-posedness of CSS for $(\bm\phi_0, \bA_0) \in H^{\frac{3}{4}+ } \times H^{\frac{1}{4} + }$, and the global solutions was also established for $(\bm\phi_0, \bA_0) \in H^{1} \times H^{\frac{1}{2}}$ by energy conservation. However, in 2-dimensional case, the well-posed problem is more difficult, and there is no well-posedness results on \eqref{E-L}. On the other hand, low regularity problem is a key topic in the study of field equations, which has been addressed in wave maps and Chern-Simons-Higgs.

Motivated by these insightful works, we study the nonlinear structure and local well-posedness of \eqref{E-L} with low regularity data. Precisely, we prove the local well-posedness of solutions if the initial data $(\bm{\phi}_0,\bA_0)\in H^{1+}(\R^2)\times H^{\frac12+}(\R^2)$, where the critical regularity of $\bm{\phi}_0$ and $\bA_0$ is $s_c=1$ and $s'_c=0$. So our result is almost critical for the matter field. Our proof is based on identifying null forms within the system and utilizing bilinear estimates in wave-Sobolev space. We state our result as follows. 


\subsection{Statement of results}
\begin{thm}\label{mainthm}
		Let $s>1$, consider the Cauchy problem of
		Chern-Simons gauged $O(3)$ sigma equations \eqref{f0}-\eqref{lorenz}, with the initial data in the following Sobolev space:
		\begin{equation*}
			\begin{split}
				& \bA(0,\cdot)=\ba\in H^{s-\frac12}(\R^2),\quad \partial_t \bA (0,\cdot)=\ba'\in H^{s-\frac32}(\R^2),
				\\
				& \bm{\phi}(0,\cdot)=\bm{\phi}_0\in H^s(\R^2),\quad \p_t\bm{\phi}(0,\cdot)=\bm{\phi}_1\in H^{s-1}(\R^2)
			\end{split}
		\end{equation*}
		satisfying the constraint \eqref{constraint1} and $\langle\bm{\phi}_0,\,\bm{\phi}_1\rangle=0$, where
		\begin{equation*}
			\ba=(a_0,a_1,a_2), \quad \ba'=(a'_0,a'_1,a'_2). 
		\end{equation*}
		Then there exists a unique solution $\bA$ and $\bm{\phi}$ to \eqref{f0}-\eqref{f3} on $[0,T]\times \R^2$ $(T>0$ and $T$ depends on $\|\ba\|_{H^{s-\frac12}},\ \|\ba'\|_{H^{s-\frac32}},\|\bm{\phi}_0\|_{H^s},\,\|\bm{\phi}_1\|_{s-1}$, and $m)$. under the Lorenz gauge condition $\partial_\mu A^\mu=0$ satisfying 
		\begin{equation*}
			\begin{split}
				& \bA\in C([0,T];H^{s-\frac12}(\R^2))\cap C^1([0,T];H^{s-\frac32}(\R^2)),
				\\
				& \bm{\phi}\in C([0,T];H^s(\R^2))\cap C^1([0,T];H^{s-1}(\R^2)).
			\end{split}
		\end{equation*}
\end{thm}

\begin{remark}
    The proof of Theorem \ref{mainthm} will use a standard approximation procedure, other standard statements of local well-posedness (such as uniqueness in the iteration space, continuous dependence on the data) follow easily. We omit the precise statements for brevity.

    \begin{remark}
    	For the critical regularity of matter fields is $s_c=1$, so our Theorem \ref{mainthm} is almost critical for the matter fields. Whether the result of Theorem \ref{mainthm} is sharp or not, we could not answer it here. We expect that the local well-posedness of Chern-Simons gauged $O(3)$ sigma equations also holds if $(\bm{\phi}_0,\bA_0) \in H^{1+} \times H^{0+}$.
    \end{remark}

\end{remark}
\begin{remark}
One cannot use $H^s$ spaces for functions with target $\mathcal{S}^2$, since they do not belong to $L^2$. Thus when we say that $\bm{\phi}(t,\,\cdot)\in H^s$, we mean that  $\bm{\phi}(t,\,\cdot)-n_3 \in H^s$.
\end{remark}


\subsection{Notations}
We first give the definition of space $H^{s,\,b}$ as well as we used. We define the function spaces with the Fourier transforms on $\mathbb{R}^{1+n}$, which denoted by $\widehat{\cdot}$ and $\widetilde{\cdot}$:
\begin{equation*}
	\widehat{f}(\xi)=\int_{\mathbb{R}^n} \text{e}^{\text{i}x\cdot \xi}f(x)dx, \quad 	\widetilde{F}(\tau,\xi)=\int_\mathbb{R} \int_{\mathbb{R}^n} \text{e}^{\text{i}(t\tau+x\cdot \xi)}F(t,x)dx dt.
\end{equation*}

For $s, b \in \mathbb{R}$, let $H^{s,b}$ be the space\footnote{The space $\mathcal{S}'(\mathbb{R}^{1+n})$ is the dual space of Schwartz functions.}
\begin{equation*}
	H^{s,b}=\left\{ u\in \mathcal{S}'(\mathbb{R}^{1+n}): \left< \xi \right>^s \left<||\tau|-|\xi||\right>^{b} \widetilde{u}(\tau,\xi) \in L^2(\mathbb{R}^{1+n}) \right\},
\end{equation*}
where we set $\left< \xi \right>=1+|\xi|$ and  $\left<||\tau|-|\xi||\right>=1+ ||\tau|-|\xi||$. For $b>\frac12$, then $H^{s,b} \hookrightarrow H^s$(c.f. \cite{Selberg}). We use the notation $\|f\|_{s,b}$ to denote a norm in $H^{s,b}$, that is
\begin{equation}\label{HS1}
	\|u\|_{s,b}=\| \left< \xi \right>^s \left<||\tau|-|\xi||\right>^{b} \tilde{u}(\tau,\xi) \|_{ L^2(\mathbb{R}^{1+n}) }.
\end{equation}
We also introduce a norm
\begin{equation}\label{HS2}
	|f|_{s,b}=\|f\|_{s,b}+ \|\partial_t f\|_{s-1,b}.
\end{equation}

The operators $\Lambda$, $D$, $D_{\pm}$, and $\Lambda_{\pm}$ are denoted by
\begin{equation*}
	\begin{split}
		&\widehat{\Lambda^\alpha f}(\xi)=\left< \xi \right>^\alpha \hat{f}(\xi),\qquad \qquad \qquad \widehat{D^\alpha f}(\xi)=|\xi|^\alpha \hat{f}(\xi),
		\\
		&\widehat{D_{\pm}^\alpha F}(\tau,\xi)=||\tau|\pm|\xi||^\alpha \tilde{F}(\tau,\xi), \quad \widehat{\Lambda_{\pm}^\alpha F}(\tau,\xi)=\left<||\tau|\pm|\xi||\right>^\alpha \tilde{F}(\tau,\xi),
	\end{split}
\end{equation*}
The space $H^s(\mathbb{R}^n)$ is defined by
\begin{equation*}
	H^{s}(\mathbb{R}^n)=\left\{ u\in \mathcal{S}'(\mathbb{R}^{n}): \left< \xi \right>^s  \hat{u}(\xi) \in L^2(\mathbb{R}^{n}) \right\}.
\end{equation*}
We use $a\lesssim b$ denote $a\leq Cb$ for some constant $C$. A point in the $2+1$ dimensional Minkowski space is written as
$(t,x)=(x^\alpha)_{0\leq\alpha\leq2}$. Greek indices range from 0 to 2, and Roman indices range from 1 to 2. We raise and lower indices
with the Minkowski metric, diag$(1,-1,-1)$. We write $\partial_\alpha=\partial_{x^\alpha}$ and $\partial_t=\partial_0,$ and we also use
the Einstein notation. Therefore, $\partial^i\partial_i=\Delta$, and $\partial^\alpha\partial_\alpha=\partial_t^2-\Delta=\Box$. We also set $\partial=(\partial_t, \partial_1, \partial_2)$.

\subsection{Organization of the rest paper} In next Section \ref{sec:2}, we introduce bilinear estimates and null structure of the \eqref{E-L} model. In Section \ref{sec:3}, we provide a proof of Theorem\ref{mainthm}. In the Appendix \ref{sec:4}, we derive the energy conservation of the \eqref{E-L}.

\section{Preliminaries}\label{sec:2}
In this part, we introduce some known bilinear estimates and discuss null structure, which plays a crucial role for proving the main theorem.
\subsection{Bilinear estimates}

\ \\
In this subsection, we introduce several useful lemmas, including product estimates of the form  $H^{s_1,b_1}\cdot H^{s_2,b_2}\hookrightarrow H^{-s_0,-b_0},$ which imply
$$\|uv\|_{H^{-s_0,-b_0}}\leq C\|u\|_{H^{s_1,b_1}}\|v\|_{H^{s_2,b_2}}\,\,\mbox{for all}\,\,u,v\in\mathcal{S}(\R^{1+n})$$
where $C$ depends on the $s_\alpha,b_\alpha$ and $d$. When this inequality holds, we say that the exponent matrix
$$
\begin{pmatrix}
s_0&s_1&s_2\\
b_0&b_1&b_2\\
\end{pmatrix}
$$
is a product. In a recent paper \cite{DFS}, the following product estimate in $\R^{1+2}$ is established.
\begin{lem}[\cite{DFS}, Theorem 5.1]\label{product}
Assume
\begin{align*}
& b_0+b_1+b_2>\frac12,  \\
&  b_0+b_1\geq 0,\\
& b_1+b_2\geq 0, \\
& b_0+b_2\geq 0,\\
  & s_0+s_1+s_2>\frac32-(b_0+b_1+b_2),\\
  & s_0+s_1+s_2>1-\mathrm{min}(b_0+b_1,b_0+b_2,b_1+b_2),\\
  & s_0+s_1+s_2>\frac12-\mathrm{min}(b_0,b_1,b_2),\\
  & s_0+s_1+s_2>\frac34,\\
     & (s_0+b_0)+2s_1+2s_2>1,\\
     & 2s_0+(s_1+b_1)+2s_2>1,\\
     & 2s_0+2s_1+(s_2+b_2)>1,\\
        & s_0+s_1\geq \mathrm{max} (0,\, -b_2),\\
        & s_1+s_2\geq \mathrm{max}(0,\, -b_0),\\
        & s_0+s_2\geq \mathrm{max} (0,\, -b_1).
\end{align*}
Then
$$
\begin{pmatrix}
s_0&s_1&s_2\\
b_0&b_1&b_2\\
\end{pmatrix}
$$
is a product.
\end{lem}

\begin{lem}[\cite{Selberg},Theorem 12, Theorem 13]\label{nonlinearE}
	Assume $s\in \mathbb{R}$, $\theta \in (\frac12,1)$, $\epsilon \in [0,1-\theta]$. Consider the Cauchy problem for the linear wave equation
	\begin{equation}\label{linearw}
		\begin{cases}
			& \square u=F, \quad (t,x)\in \mathbb{R}^{1+n},
			\\
			& u|_{t=0} =f, \quad  \partial_t u|_{t=0} = g,
		\end{cases}
	\end{equation}
where $f,g$ and $F$ satisfy $f\in H^s$, $g\in H^{s-1}$, and $F\in H^{s-1,\theta+\epsilon-1}$. 
Define a cut-off function $\chi$ satisfying
\[
	\chi \in C^\infty_c(\mathbb{R}), \quad \chi=1 \ \text{on} \ [-1,1], \quad \text{supp} \chi \subseteq (-2,2),
\]
and another cut-off function $\eta$ such that
\[
	\eta \in C^\infty_c(\mathbb{R}), \quad \eta=1 \ \text{on} \ [-2,2], \quad \text{supp} \ \eta \subseteq (-4,4).
\]
 
 For $0<T<1$, define $u(t)$ by
	\begin{equation}\label{defu}
		u(t)=\chi(t)u_0+ \chi(\frac{t}{T})u_1+u_2,
	\end{equation}
	where
	\begin{equation}\label{defu0}
		\begin{split}
			u_0= &\cos(tD)f +  D^{-1} \sin(tD)g,
			\\
			F_1= & \eta( T^{\frac12}\Lambda_{-} )F, \quad F_2= ( 1-\eta( T^{\frac12}\Lambda_{-} ) )F,
			\\
			u_1= & \int^t_0  D^{-1} \sin( (t-t')D ) F_1(t')dt',
			\\
			u_2= & \square^{-1} F_2.
		\end{split}
	\end{equation}
	Then, the function $u$ defined in \eqref{defu}-\eqref{defu0} satisfies the following estimate
	\begin{equation}\label{none1}
		|u|_{s,\theta} \leq C_0(\|f\|_{H^s}+ \|g\|_{H^{s-1}}+ T^{\frac{\epsilon}{2}}\| F\|_{s-1,\theta+\epsilon-1} ),
	\end{equation}
	where $C_0$ only depends on $\chi$ and $\theta$. Moreover, $u$ is the unique solution of \eqref{linearw} on $[0,T]\times \mathbb{R}^n$ such that $u\in C([0,T];H^s) \cap C^1([0,T];H^{s-1})$.
\end{lem}
\begin{remark}
	For $\epsilon=0$ or $\epsilon\in (0,1-\theta]$, we refer the readers to Selberg's paper \cite{Selberg} Theorem 12 and Theorem 13 respectively. 
\end{remark}
Let $u$ and $v$ be scalar functions satisfying
\begin{equation}\label{eu}
	\begin{cases}
		& \square u=0, \quad (t,x)\in \mathbb{R}^{1+n},
		\\
		& u|_{t=0} =u_0, \quad  \partial_t u|_{t=0} = u_1,
	\end{cases}
\end{equation}
and
\begin{equation}\label{ev}
	\begin{cases}
		& \square v=0, \quad (t,x)\in \mathbb{R}^{1+n},
		\\
		& v|_{t=0} =v_0, \quad  \partial_t v|_{t=0} = v_1.
	\end{cases}
\end{equation}

For simplicity, we take $u_1=v_1=0$ in \eqref{eu}-\eqref{ev}. The following null form estimates in Sobolev space which was proved by Foschi and Klainerman in $n$-dimensions ($n\geq 2$) in \cite{FK}.
\begin{lem}[\cite{FK}, Corollary 13.3]\label{nq0}
The estimate
\begin{equation}\label{nq00}
	\| D^{\beta_0} D_{+}^{\beta_{+}} D_{-}^{\beta_{-}} Q_0(u,v) \|_{L^2(\mathbb{R}^{1+n})} 
	\lesssim \| D^{\alpha_1} f_1 \|_{L^2(\mathbb{R}^n)} \| D^{\alpha_2} f_2 \|_{L^2(\mathbb{R}^n)},
\end{equation}
holds when $\alpha_1, \alpha_2, \beta_0, \beta_{+}$, and $\beta_{-}$ satisfy the following conditions:
\begin{align*}
& \beta_0+\beta_{+}+\beta_{-}=\alpha_1+\alpha_2-\frac{n+3}{2},
\\
& \beta_{-}\geq -\frac{n+1}{4}, \quad \beta_{0}> -\frac{n-1}{2},
\\
& \alpha_1+\alpha_2\geq \frac12, \quad \alpha_i \leq \beta_{-} + \frac{n+1}{2}, , \quad i=1,2
\\
& (\alpha_1+\alpha_2,\beta_{-} ) \neq (\frac12,- \frac{n+1}{4}), \quad (\alpha_i,\beta_{-})\neq (\frac{n+1}{4}, - \frac{n+1}{4}).
\end{align*}
\end{lem}

\begin{lem}[\cite{FK}, Corollary 13.4]\label{nq1}
	The estimate
	\begin{equation}\label{nq01}
		\| D^{\beta_0} D_{+}^{\beta_{+}} D_{-}^{\beta_{-}} Q_{ij}(u,v) \|_{L^2(\mathbb{R}^{1+n})} 
		\lesssim \| D^{\alpha_1} f_1 \|_{L^2(\mathbb{R}^n)} \| D^{\alpha_2} f_2 \|_{L^2(\mathbb{R}^n)},
	\end{equation}
	holds when $\alpha_1, \alpha_2, \beta_0, \beta_{+}$, and $\beta_{-}$ satisfy the following conditions:
	\begin{align*}
		& \beta_0+\beta_{+}+\beta_{-}=\alpha_1+\alpha_2-\frac{n+3}{2},
		\\
		& \beta_{-}\geq -\frac{n-1}{4}, \quad \beta_{0}> -\frac{n+1}{2},
		\\
		& \alpha_1+\alpha_2\geq \frac32, \quad \alpha_i \leq \beta_{-} + \frac{n+1}{2}, \quad i=1,2,
		\\
		& (\alpha_1+\alpha_2,\beta_{-} ) \neq (\frac12,- \frac{n-1}{4}), \quad (\alpha_i,\beta_{-})\neq (\frac{n+3}{4}, - \frac{n-1}{4}).
	\end{align*}
\end{lem}

\begin{lem}[\cite{FK}, Corollary 13.5]\label{nq2}
	The estimate
	\begin{equation}\label{nq02}
		\| D^{\beta_0} D_{+}^{\beta_{+}} D_{-}^{\beta_{-}} Q_{0j}(u,v) \|_{L^2(\mathbb{R}^{1+n})} 
		\lesssim \| D^{\alpha_1} f_1 \|_{L^2(\mathbb{R}^n)} \| D^{\alpha_2} f_2 \|_{L^2(\mathbb{R}^n)},
	\end{equation}
	holds when $\alpha_1, \alpha_2, \beta_0, \beta_{+}$, and $\beta_{-}$ satisfy the following conditions:
	\begin{align*}
		& \beta_0+\beta_{+}+\beta_{-}=\alpha_1+\alpha_2-\frac{n+3}{2},
		\\
		& \beta_{-}\geq -\frac{n+1}{4}, \quad \beta_{0}> -\frac{n-1}{2},
		\\
		& \alpha_1+\alpha_2\geq \frac12, \quad \alpha_i \leq \beta_{-} + \frac{n+1}{2},
		\\
		& (\alpha_1+\alpha_2,\beta_{-} ) \neq (\frac12,- \frac{n+1}{4}), \quad (\alpha_i,\beta_{-})\neq (\frac{n+1}{4}, - \frac{n+1}{4}).
	\end{align*}
\end{lem}

\subsection{Null structure}

\ \\
In this subsection, we introduce the null structure of $A^\mu \partial_\mu  \bm{\phi}$. We define the Riesz transform $R_i=D^{-1}\partial_i$ ($i=1,2$) and $\mathring{\bA}=(A_1,A_2)$. Following Pecher's work \cite{Pecher} on Maxwell-Klein-Gordon system, we can decompose the vector $A$ by two parts: divergence free and curl free. Thus, we have
\begin{equation*}
	\mathring{\bA}=\mathring{\bA}^{df}+\mathring{\bA}^{cf},
\end{equation*}
where $\mathring{\bA}^{df}=(A_1^{df},A_2^{df})$ and $\mathring{\bA}^{cf}=(A_1^{cf},A_2^{cf})$, and
\begin{equation*}
	A_i^{df}=R^k( R_i A_k - R_k A_i), \quad A_j^{cf}=-R_i R^k A_k, \quad i=1,2.
\end{equation*}
Therefore we obtain
\begin{equation}\label{dec}
	A^\mu \partial_\mu \bm{\phi}=\underbrace{(-A^0 \partial_t \bm{\phi} + \mathring{\bA}^{cf} \cdot \nabla \bm{\phi})}_{\equiv P_1}+ \underbrace{\mathring{\bA}^{df} \cdot \nabla \bm{\phi}}_{\equiv P_2}.
\end{equation}
In the Lorentz gauge $\partial_\mu A^\mu=0$, we can compute out 
\begin{equation}\label{P1}
	\begin{split}
		P_1=& -A_0 \partial_t \bm{\phi} - R_i R_k A^k \partial^i \bm{\phi}
		\\
		=& -A_0 \partial_t \bm{\phi}- ( D^{-2} \nabla \partial_t A_0 ) \cdot \nabla \bm{\phi}
		\\
		=& \partial_i (D^{-1}R^i  A_0 ) \partial_t \bm{\phi} - \partial_t ( D^{-1}R_i A_0 ) \partial^i \bm{\phi}
		\\
		=& Q_{0i} (D^{-1} R^i A_0,\bm{\phi} ).
	\end{split}
\end{equation}
Similarly, we also have
\begin{equation}\label{P2}
	\begin{split}
		P_2=& \mathring{\bA}^{df} \cdot \nabla \bm{\phi}
		\\
		=& R^k( R_i A_k- R_k A_i ) \partial^i \bm{\phi}
		\\
		=& D^{-2} \partial^k \partial_i A_k \partial^i \bm{\phi} + A_i \partial^i \bm{\phi}
		\\
		=& -\frac12 \left(  D^{-2}  (\partial_i \partial^i A_k - \partial_i \partial_k A^i ) \partial^k \bm{\phi} -D^{-2} (\partial^k \partial_i A_k - \partial_k \partial^k A_i ) \partial^i \bm{\phi}  \right)
		\\
		=& -\frac12 \left(  \partial_i  D^{-1}  (R^i A_k - R_k A^i ) \partial^k \bm{\phi} - \partial^k  D^{-1}  (R_i A_k - R_k A_i ) \partial^i \bm{\phi}  \right)
		\\
		=& -\frac12 Q_{ik} (D^{-1} (R^i A^k - R^k A^i),\bm{\phi} ).
	\end{split}
\end{equation}
Adding \eqref{P1} and \eqref{P2} to \eqref{dec}, we get
\begin{equation}\label{decf}
	A^\mu \partial_\mu \bm{\phi}= -\frac12 Q_{ik} (D^{-1} (R^i A^k - R^k A^i),\bm{\phi} )+ Q_{0i} (D^{-1} R^i A_0,\bm{\phi} ).
\end{equation}


\section{Proof of Theorem \ref{mainthm}}\label{sec:3}
This section is devoted to the proof of Theorem \ref{mainthm}. Since \eqref{c}, so the specific value of the positive constant $\kappa$ does not affect our analysis, so we assume $\kappa=1$ for simplicity.
\subsection{A Prior Energy Estimates.}\label{s1}
Recall \eqref{main1}-\eqref{main2}. To be simple, we set
\begin{equation*}
	M_0=\|\bm{\phi}_0\|_{H^{s}} + \|\bm{\phi}_1\|_{H^{s-1}} +\|\ba\|_{H^{s-\frac12}} + \|\ba'\|_{H^{s-\frac32}},
\end{equation*}
and
\begin{equation}\label{GH}
	\begin{split}
		\bG =&-\phi\left(
		Q_0(\bm{\phi},\bm{\phi})+2A_{\mu}\p^{\mu}\bm{\phi}\cdot(n_3\times\bm{\phi})+A_{\mu} A^{\mu}|n_3\times\bm{\phi}|^2\right) -2 A^{\mu} (n_3\times \p_{\mu} \bm{\phi}) 
		\\
		& -A_\mu A_\mu n_3\times(n_3\times\bm{\phi})-\frac{1}{\kp^2}(\bm{\phi}(n_3\cdot\bm{\phi})-n_3(\bm{\phi}\cdot\bm{\phi}))(1-n_3\cdot\bm{\phi})^2(1+2n_3\cdot\bm{\phi}),
		\\
		\bH=& (H_0,H_1,H_2), \quad H_\mu=\frac{1}{\kappa} \epsilon_{\mu\nu\rho}\left( Q^{\nu \rho}(\bm{\phi}, \bm{\phi} \times n_3) + A^{\rho} \p^{\nu} \left( (\bm{\phi} \times (n_3 \times \bm{\phi})) \cdot n_3 \right) \right).
	\end{split}
\end{equation}		
For $s>1$ and $b=\frac{s}{2}$, assume the solution $(\bm{\phi}, \bA)$ of \eqref{main1}-\eqref{main2} satisfying
\begin{equation}\label{Xs}
	|\bm{\phi}|_{s,b}+|\bA|_{s-\frac12,b}  \leq 2(1+C_0)M_0.
\end{equation}
Since Lemma \ref{nonlinearE} and Duhamel's principle, it yields
\begin{equation}\label{e00}
	\begin{split}
		& |\bm{\phi}|_{s,b}  \leq C( \|\bm{\phi}_0\|_{H^{s}} + \|\bm{\phi}_1\|_{H^{s-1}} 
		+ T^{\frac{\epsilon}{2}} \| \bG\|_{s-1,b-1+\epsilon} ),
		\\
		&  |\bA|_{s-\frac12,b}  \leq C( \| \ba \|_{H^{s-\frac12}}  +  \|\ba'\|_{H^{s-\frac32}}
		+ T^{\frac{\epsilon}{2}} \| \bH \|_{s-\frac32,b-1+\epsilon} ).
	\end{split}
\end{equation}
Next, let us bound $\| \bG\|_{s-1,b-1+\epsilon}$ and $\| \bH \|_{s-\frac32,b-1+\epsilon}$. By Lemma \ref{product}, we have
\begin{equation}\label{e01}
	\begin{split}
		\| \phi A_{\mu} A^{\mu}|n_3\times\bm{\phi}|^2 \|_{s-1,b-1+\epsilon} 
		\leq  C  | \bm{\phi} |^2_{s,b}  | \bA |^2_{s-\frac12,b},
	\end{split}
\end{equation}
and
\begin{equation}\label{e02}
	\begin{split}
		\| A_\mu A^\mu n_3\times(n_3\times\bm{\phi})\|_{s-1,b-1+\epsilon} \leq  C  | \phi |_{s,b}  | \bA |^2_{s-\frac12,b}.
	\end{split}
\end{equation}
and
\begin{equation}\label{e04}
	\begin{split}
		\|\frac{1}{\kp^2}(\bm{\phi}(n_3\cdot\bm{\phi})-n_3(\bm{\phi}\cdot\bm{\phi}))(1-n_3\cdot\bm{\phi})^2(1+2n_3\cdot\bm{\phi}) \|_{s-1,b-1+\epsilon}\leq C  (| \bm{\phi} |^2_{s,b}+ |\bm{\phi} |^5_{s,b}),
	\end{split}
\end{equation}
and
\begin{equation}\label{e06}
	\begin{split}
		\|\bA \partial \bm{\phi} \|_{s-\frac32,b-1+\epsilon}\leq C | \bm{\phi} |_{s,b}  | \bA |^2_{s-\frac12,b} . 
	\end{split}
\end{equation}
By Lemma \ref{product} and Lemma \ref{nq0}(taking $\beta_{-}=b-1+\epsilon, \beta_{+}=0, \beta_0=s-1, \alpha_1=s,\alpha_2=\frac12+b-\epsilon$), we can derive that
\begin{equation}\label{e05}
	\begin{split}
		\| \bm{\phi} Q_0(\bm{\phi},\bm{\phi})\|_{s-1,b-1+\epsilon} \leq & C\| \bm{\phi} \|_{s,b} \|Q_0(\bm{\phi},\bm{\phi})\|_{s-1,b-1+\epsilon}
		\\
		\leq & C\| \bm{\phi} \|_{s,b} \cdot \| \bm{\phi} \|_{s,b} \| \bm{\phi} \|_{\frac12+b-\epsilon,b}
		\\
		\leq & C| \bm{\phi} |^3_{s,b},
	\end{split}
\end{equation}
for $b=\frac{s}{2}<s-\frac12$.

Now our goal transfers to $A^\mu \partial_\mu \bm{\phi}$ and $Q(\bm{\phi},\bm{\phi})$. From \eqref{decf}, there is a very good null structure for $A^\mu \partial_\mu \bm{\phi}$. But there is an operator $D^{-1}$ in \eqref{decf}. Roughly speaking, it may cause singularity for low frequency. So we decompose it into low and high frequency. Set $\Omega=\{ \xi \in \mathbb{R}^2: |\xi| \leq 1 \}$. By \eqref{dec}, we have 
\begin{equation}\label{e07}
	\begin{split}
		\| A^\mu \partial_\mu \bm{\phi}\|_{s-1,b-1+\epsilon} \leq & \| P_1\|_{s-1,0}+ \| P_2\|_{s-1,0}
		\\
		=& \| \Lambda^{s-1} P_1\|_{L^2(\mathbb{R} \times \Omega)}+ \| \Lambda^{s-1} P_1\|_{L^2(\mathbb{R} \times (\mathbb{R}^2\backslash \Omega))}
		\\
		& + \| \Lambda^{s-1} P_2\|_{L^2(\mathbb{R} \times \Omega)}+ \| \Lambda^{s-1} P_2\|_{L^2(\mathbb{R} \times (\mathbb{R}^2\backslash \Omega))}.
	\end{split}
\end{equation}
Seeing \eqref{P1}-\eqref{P2}, by Bernstein's inequality and product estimates, we obtain
\begin{equation}\label{e09}
	\begin{split}
		\| \Lambda^{s-1} P_1\|_{L^2(\mathbb{R} \times \Omega)} + \| \Lambda^{s-1} P_2\|_{L^2(\mathbb{R} \times \Omega)} \leq  & C\|\Lambda^{-\frac12} \left(  -A_0 \partial_t \bm{\phi} - R_i R_k A^k \partial^i \bm{\phi} \right)\|_{L^2(\mathbb{R} \times \Omega)}
		\\
		& + \| \Lambda^{-\frac12} \left( R^k( R_i A_k- R_k A_i ) \partial^i \bm{\phi} \right)\|_{L^2(\mathbb{R} \times \Omega)}
		\\
		\leq &  C\|\Lambda^{\frac12} \bA \|_{L^2(\mathbb{R} \times \Omega)}\|\Lambda^{s-1} \partial_t \bm{\phi} \|_{L^2(\mathbb{R} \times \Omega)}
		\\
		\leq &  C | \bA |_{s-\frac12,b} | \bm{\phi} |_{s,b}.
	\end{split}
\end{equation}
Using Bernstein inequality and \eqref{P1}-\eqref{P2} again, we can show that
\begin{equation}\label{e06a}
	\begin{split}
		& \| \Lambda^{s-1} P_1\|_{L^2(\mathbb{R} \times (\mathbb{R}^2\backslash \Omega))}+ \| \Lambda^{s-1} P_2\|_{L^2(\mathbb{R} \times (\mathbb{R}^2\backslash \Omega))}
		\\
		\leq & C\| D^{s-1} Q_{ik} \big(D^{-1} (R^i A^k - R^k A^i),\bm{\phi} \big)\|_{0,0}+ \| D^{s-1} Q_{0i} \left(D^{-1} R^i A_0,\bm{\phi} \right) \|_{0,0}.
	\end{split}
\end{equation}
Due to Lemma \ref{nq2} (taking $\beta_{-}=0, \beta_{+}=0, \beta_0=s-1, \alpha_1=\frac32, \alpha_2=s$), it yields
\begin{equation}\label{e08}
	\begin{split}
		& \| D^{s-1} Q_{ik} \big(D^{-1} (R^i A^k - R^k A^i),\bm{\phi} \big)\|_{0,0}+ \| D^{s-1} Q_{0i} \left(D^{-1} R^i A_0,\bm{\phi} \right) \|_{0,0}
		\\
		\leq & C\| D^{\frac32} (D^{-1} \bA)\|_{0,b} \| D^{s}\bm{\phi} \|_{0,b}
		\\
		\leq & C | \bA |_{s-\frac12,b} | \bm{\phi} |_{s,b}.
	\end{split}
\end{equation}
Combining with \eqref{e06a} and \eqref{e08}, we can see
\begin{equation}\label{e11}
	\begin{split}
		\| \Lambda^{s-1} P_1\|_{L^2(\mathbb{R} \times (\mathbb{R}^2\backslash \Omega))}+ \| \Lambda^{s-1} P_2\|_{L^2(\mathbb{R} \times (\mathbb{R}^2\backslash \Omega))}
		\leq  C | \bA |_{s-\frac12,b} | \bm{\phi} |_{s,b}.
	\end{split}
\end{equation}
Adding \eqref{e09} and \eqref{e11} to \eqref{e07}, we have proved
\begin{equation}\label{a12}
	\| A^\mu \partial_\mu \bm{\phi}\|_{s-1,b-1+\epsilon} \leq C | \bA |_{s-\frac12,b} | \bm{\phi} |_{s,b}.
\end{equation}
For $Q(\bm{\phi},\bm{\phi})$, using Lemma \ref{nq1}(taking $\beta_{-}=0, \beta_{+}=0, \beta_0=s-\frac32, \alpha_1=1, \alpha_2=s$), we can infer from \eqref{decf} that
\begin{equation}\label{e15}
	\begin{split}
		\| Q(\bm{\phi}, \bm{\phi})\|_{s-\frac32,b-1+\epsilon} \leq & \| Q( \bm{\phi}, \bm{\phi})\|_{s-\frac32,0}
		\\
		\leq & C|\bm{\phi} |_{s,b} |\bm{\phi} |_{1,b} 
		\\
		\leq & C|\bm{\phi} |^2_{s,b}.
	\end{split}
\end{equation}
Summing up our outcome \eqref{e01}, \eqref{e02}, \eqref{e04}, \eqref{e06}, \eqref{e05}, \eqref{a12}, and \eqref{e15}, we can conclude that
\begin{equation}\label{e17}
	\begin{split}
		& \| \bG \|_{s-1,b-1-\epsilon} + \| \bH \|_{s-\frac32,b-1-\epsilon} 
		\\
		\leq & C ( |\bm{\phi} |^2_{s,b}|\bA |^2_{s-\frac12,b}+ |\bm{\phi} |_{s,b}|\bA |^2_{s-\frac12,b}+|\bm{\phi} |^3_{s,b}+ |\bA |^2_{s-\frac12,b}+ |\bm{\phi} |^2_{s,b} ).
	\end{split}
\end{equation}
By \eqref{e00} and \eqref{e17}, 
\begin{equation*}\label{T}
	\begin{split}
		& |\bm{\phi} |_{s,b}+ |\bA |_{s-\frac12,b} 
		\\
		\leq & C_0 M_0+ C_0 C ( |\bm{\phi} |^2_{s,b}|\bA |^2_{s-\frac12,b}+ |\bm{\phi} |_{s,b}|\bA |^2_{s-\frac12,b}+|\bm{\phi} |^3_{s,b}+ |\bA |^2_{s-\frac12,b}+ |\bm{\phi} |^2_{s,b} ).
	\end{split}	
\end{equation*}
Due to the standard continuation argument, the uniform bound $|\bm{\phi} |_{s,b}+ |\bA |_{s-\frac12,b} \leq (1+C_0) M_0$ follows when the positive number $T$ is sufficiently small. Moreover, the uniqueness of solution follows
from a routine argument and the energy estimates derive above. We omit the details for brevity.
\subsection{Local Existence with finite energy.}
After obtaining the uniform energy bounds, we now give a proof of local existence for \eqref{main1}-\eqref{main2}. For \eqref{main1}-\eqref{main2} is a semi-linear system, we will obtain the local existence through a standard
approximation procedure. Firstly, we define
\begin{equation}\label{e20}
	\begin{cases}
		& \Box \bm{\phi}^{(1)}=0,  \quad (t,x) \in [0,T] \times \mathbb{R}^2,
		\\
		& \bm{\phi}^{(1)}|_{t=0}=  \bm{\phi}_0, \quad  \partial_t  \bm{\phi}^{(1)} |_{t=0}= \bm{\phi}_1,
	\end{cases}
\end{equation}
and
\begin{equation}\label{e21}
	\begin{cases}
		& \Box \bA^{(1)}=0,  \quad (t,x) \in [0,T] \times \mathbb{R}^2,
		\\
		& \bA^{(1)}|_{t=0}=  \ba, \quad  \partial_t  \bA^{(1)} |_{t=0}= \ba',
	\end{cases}
\end{equation}
For each positive integer $m\geq 1$, we define
\begin{equation}\label{e22}
	\begin{cases}
		& \Box \bm{\phi}^{(m+1)}=\bG^{(m)},  \quad (t,x) \in [0,T] \times \mathbb{R}^2,
		\\
		& \bm{\phi}^{(m+1)}|_{t=0}=  \bm{\phi}_0, \quad  \partial_t  \bm{\phi}^{(m+1)} |_{t=0}= \bm{\phi}_1,
	\end{cases}
\end{equation}
and
\begin{equation}\label{e24}
	\begin{cases}
		& \Box A_\mu^{(m+1)}=H_\mu^{(m)},  \quad (t,x) \in [0,T] \times \mathbb{R}^2,
		\\
		& A_\mu^{(m+1)}|_{t=0}=  a_\mu, \quad  \partial_t  A_\mu^{(m+1)} |_{t=0}= a'_\mu,
	\end{cases}
\end{equation}
where $\bG^{(m)}$ and $H^{(m)}_\mu$ have the same formulation with $\bG$ and $H_\mu$ only by replacing $(\bm{\phi},\bA)$ to $(\bm{\phi}^{(m)},\bA^{(m)})$. We also set the three component vector $\bH^{(m)}= (H^{(m)}_0,H^{(m)}_1,H^{(m)}_2)$.

By Lemma \ref{nonlinearE}, then the problem \eqref{e20}, \eqref{e21}, \eqref{e22}, and \eqref{e24} all have a unique solution. Next, we will prove the iteration sequence $\{(\bm{\phi}^{(m)}, \bA^{(m)})\}_{m\geq 1}$ is convergent. By \eqref{e22} and \eqref{e24}, then the minus $(\bm{\phi}^{(m+1)}-\bm{\phi}^{(m)}, \bA^{(m+1)}-\bA^{(m)})$ satisfies
\begin{equation}\label{e25}
	\begin{cases}
		& \Box (\bm{\phi}^{(m+1)}-\bm{\phi}^{(m)})=\bG^{(m)}-\bG^{(m-1)},  \quad (t,x) \in [0,T] \times \mathbb{R}^2,
		\\
		& \bm{\phi}^{(m+1)}-\bm{\phi}^{(m)}|_{t=0}= 0, \quad  \partial_t  (\bm{\phi}^{(m+1)}-\bm{\phi}^{(m)}) |_{t=0}= 0,
	\end{cases}
\end{equation}
and
\begin{equation}\label{e26}
	\begin{cases}
		& \Box (\bA^{(m+1)}-\bA^{(m)})=\bH^{(m)}-\bH^{(m-1)},  \quad (t,x) \in [0,T] \times \mathbb{R}^2,
		\\
		& \bA^{(m+1)}-\bA^{(m)}|_{t=0}= 0, \quad  \partial_t  (\bA^{(m+1)}-\bA^{(m)}) |_{t=0}= 0.
	\end{cases}
\end{equation}
For \eqref{e25} and \eqref{e26}, using Lemma \ref{nonlinearE} again, we have
\begin{equation*}\label{e28}
	\footnotesize	\begin{split}
		& |\bm{\phi}^{(m+1)}-\bm{\phi}^{(m)} |_{s,b}+ |\bA^{(m+1)}-\bA^{(m)} |_{s-\frac12,b} 
		\\
		\leq & C T^{\frac{\epsilon}{2}}  
		|\bm{\phi}^{(m)}-\bm{\phi}^{(m-1)}|_{s,b} ( |\bm{\phi}^{(m)} |_{s,b} |\bA^{(m)} |^2_{s-\frac12,b}+ |\bA^{(m)} |^2_{s-\frac12,b}+|\bm{\phi}^{(m)} |^2_{s,b}+  |\bm{\phi}^{(m)} |_{s,b} )
		\\
		& + C T^{\frac{\epsilon}{2}}  
		|\bm{\phi}^{(m)}-\bm{\phi}^{(m-1)} |_{s,b} ( |\bm{\phi}^{(m-1)}  |_{s,b} |\bA^{(m-1)} |^2_{s-\frac12,b}+ |\bA^{(m-1)} |^2_{s-\frac12,b}+|\bm{\phi}^{(m-1)}  |^2_{s,b}+  |\bm{\phi}^{(m-1)}  |_{s,b} )
		\\
		& + C T^{\frac{\epsilon}{2}}  
		|\bA^{(m)}-\bA^{(m-1)} |_{s-\frac12,b} ( |\bm{\phi}^{(m)} |^2_{s,b} |\bA^{(m)} |_{s-\frac12,b}+ |\bm{\phi}^{(m)} |_{s,b}|\bA^{(m)} |_{s-\frac12,b}+  |\bA^{(m)} |_{s-\frac12,b} )
		\\
		& + C T^{\frac{\epsilon}{2}}  
		|\bA^{(m)}-\bA^{(m-1)}|_{s-\frac12,b} ( |\bm{\phi}^{(m-1)} |^2_{s,b} |\bA^{(m-1)} |_{s-\frac12,b}+ |\bm{\phi}^{(m-1)} |_{s,b}|\bA^{(m-1)} |_{s-\frac12,b}+  |\bA^{(m-1)}|_{s-\frac12,b} ).
	\end{split}	
\end{equation*}
Using the uniform energy \eqref{Xs} in subsection \ref{s1} and $T$ is small enough, it tells us
\begin{equation*}\label{e29}
	\begin{split}
		|\bm{\phi}^{(m+1)}-\bm{\phi}^{(m)} |_{s,b}+ |\bA^{(m+1)}-\bA^{(m)} |_{s-\frac12,b}
		\leq \frac12
		(	|\bm{\phi}^{(m)}-\bm{\phi}^{(m-1)} |_{s,b}+ |\bA^{(m)}-\bA^{(m-1)} |_{s-\frac12,b} ).
	\end{split}	
\end{equation*}
Therefore, the sequence of solutions $\{ (\bm{\phi}^{(m)}, \bA^{(m)})  \}_{m\geq 1} $ is convergent. Let the limit be $(\bm{\phi}, \bA)$. So we have
\begin{equation*}
	\lim_{m\rightarrow \infty} | \bm{\phi}^{(m)} - \bm{\phi} |_{s,b}=0, \quad  	\lim_{m\rightarrow \infty} | \bA^{(m)} - \bA |_{s,b}=0.
\end{equation*}
This strong convergence makes it easy to check that $(\bm{\phi}, \bA)$ is the solution of \eqref{main1}-\eqref{main2} and
\begin{equation*}\label{e30}
	\begin{split}
		|\bm{\phi}|_{s,b}+ |\bA |_{s-\frac12,b}
		\leq 2(1+C_0)M_0. 
	\end{split}	
\end{equation*}
Combining the local existence and uniform bounds in Sect. \ref{s1}, we have proved Theorem \ref{mainthm}.

\section{appendix: Energy conservation}\label{sec:4}

\ \\
We present the detailed calculation of the energy conservation law as given in equation \eqref{energy}:
\[
	E(\bm{\phi}, A)=\frac12 \int_{\R^2}\left( |D_\mu\bm{\phi}|^2 +\frac{1}{\kappa^2}(1+ n_3\cdot\bm{\phi})(1-n_3\cdot\bm{\phi})^3\right)dx.
\]

We first take the inner product of both sides of $\eqref{E-L}_1$ with $D_0\bm{\phi}$ to get
\begin{align*}
    0=&D_0\bm{\phi}\cdot D_0D_0\bm{\phi}-D_0\bm{\phi}\cdot D_1D_1\bm{\phi}-D_0\bm{\phi}\cdot D_2D_2\bm{\phi}+(D_{\mu}\bm{\phi}\cdot D_{\mu}\bm{\phi})(\bm{\phi}\cdot D_{\mu}\bm{\phi})\\
    &+\frac{1}{\kappa^2}(1-n_3\cdot\bm{\phi})^2(1+2n_3\cdot\bm{\phi})\Big[(n_3\cdot\bm{\phi})(\bm{\phi}\cdot D_0\bm{\phi})-n_3\cdot D_0\bm{\phi}\big]\\
    =&\frac12\partial_0|D_0\bm{\phi}|^2-\partial_1(D_0\bm{\phi}\cdot D_1\bm{\phi})+D_1D_0\bm{\phi}\cdot D_1\bm{\phi}-\partial_2(D_0\bm{\phi}\cdot D_2\bm{\phi})+D_2D_0\bm{\phi}\cdot D_2\bm{\phi}\\
    &+|D_{\mu}\bm{\phi}|^2\partial_0\frac{|\bm{\phi}|^2}{2}+\frac{1}{\kappa^2}(1-n_3\cdot\bm{\phi})^2(1+2n_3\cdot\bm{\phi})\left[\left((n_3\cdot\bm{\phi})\partial_0\frac{|\bm{\phi}|^2}{2}\right)-\partial_0(n_3\cdot\bm{\phi})\right]\\
    =&\frac12\partial_0|D_0\bm{\phi}|^2-\partial_1(D_0\bm{\phi}\cdot D_1\bm{\phi})+D_1D_0\bm{\phi}\cdot D_1\bm{\phi}-\partial_2(D_0\bm{\phi}\cdot D_2\bm{\phi})+D_2D_0\bm{\phi}\cdot D_2\bm{\phi}\\
    &-\frac{1}{\kappa^2}(1-n_3\cdot\bm{\phi})^2(1+2n_3\cdot\bm{\phi})\partial_0(n_3\cdot\bm{\phi}),
\end{align*}
where we use the constraint condition $|\bm{\phi}|=1$. Next, we integrate the above equality over $\mathbb{R}^2$ to obtain:
\begin{align*}
    0=&\int_{\mathbb{R}^2}\frac12\partial_0|D_0\bm{\phi}|^2+D_1D_0\bm{\phi}\cdot D_1\bm{\phi}+D_2D_0\bm{\phi}\cdot D_2\bm{\phi}-\frac{1}{\kappa^2}(1-n\cdot\bm{\phi})^2(1+2n_3\cdot\bm{\phi})\partial_0(n_3\cdot\bm{\phi})dx
    \\
    =&\int_{\mathbb{R}^2}\frac12\partial_0|D_0\bm{\phi}|^2+(D_0D_1\bm{\phi}+F_{10}n_3\times\bm{\phi})\cdot D_1\bm{\phi}+(D_0D_2\bm{\phi}+F_{20}n_3\times\bm{\phi})\cdot D_2\bm{\phi}
    \\
    &+\frac12\partial_0\left(\frac{1}{\kappa^2}(1+n_3\cdot\bm{\phi})(1-n_3\cdot\bm{\phi})^3\right)dx
    \\
    =&\int_{\mathbb{R}^2}\frac12\partial_0(|D_0\bm{\phi}|^2+|D_1\bm{\phi}|^2+|D_2\bm{\phi}|^2)+F_{10}n\times\bm{\phi}\cdot D_1\bm{\phi}+F_{20}n_3\times\bm{\phi}\cdot D_1\bm{\phi}
    \\
    &+\frac12\partial_0\left(\frac{1}{\kappa^2}(1+n_3\cdot\bm{\phi})(1-n_3\cdot\bm{\phi})^3\right)dx
    \\
    =&\int_{\mathbb{R}^2}\frac12\partial_0(|D_0\bm{\phi}|^2+|D_1\bm{\phi}|^2+|D_2\bm{\phi}|^2)-\frac{1}{\kappa}(n_3\cdot\bm{\phi}\times D_2\bm{\phi})(n_3\times\bm{\phi}\cdot D_1\bm{\phi})
    \\
    &+\frac{1}{\kappa}(n_3\cdot\bm{\phi}\times D_1\bm{\phi})(n_3\times\bm{\phi}\cdot D_2\bm{\phi})
   +\frac12\partial_0\left(\frac{1}{\kappa^2}(1+n_3\cdot\bm{\phi})(1-n_3\cdot\bm{\phi})^3\right)dx
    \\
    =&\int_{\mathbb{R}^2}\frac12\partial_0(|D_0\bm{\phi}|^2+|D_1\bm{\phi}|^2+|D_2\bm{\phi}|^2)+\frac12\partial_0\left(\frac{1}{\kappa^2}(1+n_3\cdot\bm{\phi})(1-n\cdot\bm{\phi})^3\right)dx.
\end{align*}
This yields the conserved energy equation as desired.

\section*{Acknowledgement}
The work of G. Jin was supported by NNSFC under Grant Number 12201542 and the Jilin Science and Technology Development Program[YDZJ202201ZYTS311]. The work of Huali Zhang is supported by National Natural Science Foundation of China (Grant No.12101079) and the Fundamental Research Funds for the Central Universities (Grant No.531118010867).

\section*{Conflicts of interest and Data Availability Statements}
The authors declared that this work does not have any conflicts of interest. The authors also confirm that the data supporting the findings of this study are available within the article.

\end{document}